\documentclass{amsart}

\title{An Exact Value for the Ramsey Number $R(K_5, K_5-e)$}
\author{Vigleik Angeltveit}
\address{Mathematical Sciences Institute \\
Australian National University \\
Canberra, ACT 2600 \\
Australia}

\usepackage{amsxtra}
\usepackage{amsfonts}
\usepackage[latin1]{inputenc}
\usepackage{graphicx}
\usepackage{amsmath,amssymb,latexsym,amsthm,mathrsfs}
\usepackage[all]{xy}
\usepackage{pstricks}
\usepackage{verbatim}
\usepackage{hyperref}
\usepackage{placeins} 

\newtheorem{theorem}{Theorem}[section]

\theoremstyle{definition}

\makeatletter
\let\c@equation\c@theorem
\makeatother
\numberwithin{equation}{section}

 \newcommand{\cR}{\mathcal{R}}        

\begin{document}

\begin{abstract}
We compute the exact value of the Ramsey number $R(K_5, K_5-e)$. It is equal to $30$.
\end{abstract}

\maketitle

\section{Introduction}
The Ramsey number $R(s,t)$ is defined to be the smallest $n$ such that every graph of order $n$ contains either a clique $K_s$ of $s$ vertices or an independent set $\overline{K}_t$ of $t$ vertices. The Ramsey number $R(K_s, K_t-e)$ is the smallest $n$ such that every graph of order $n$ contains either a clique of $s$ vertices, or $t$ vertices with at most $1$ edge between them. See \cite{Ra94} for a survey on the currently known bounds for small Ramsey numbers.

For ease of notation, when $s$ is an integer we let $K_{s-1/2} = K_s-e$. Similarly, we extend the notation $R(s,t) = R(K_s, K_t)$ to the case when one or both of $s$ and $t$ is a half-integer. In particular we will denote $R(K_5, K_5-e)$ by $R(5,4.5)$. We also let $\cR(s,t,n)$ denote the set of isomorphism classes of Ramsey graphs of type $(s,t)$ of order $n$.

The purpose of this paper is to prove the following result:

\begin{theorem} \label{t:main2}
The Ramsey number $R(5,4.5)$ is equal to $30$.
\end{theorem}

This should be compared to the first table in \cite[Section 3.1]{Ra94}, where the value of $R(5,4.5)$ is given as being between $30$ and $33$.

The Ramsey number $R(4,5) = R(5,4) = 24$ was computed by McKay and Radziszowski \cite{McRa95} in 1995. The Ramsey number $R(5,5)$ is unknown, but is between $43$ and $46$ \cite{AnMc25}. An exact calculation of $R(5,5)$ is out of reach with current methods, and since determining $R(5,4.5)$ is between $R(5,4)$ and $R(5,5)$ in difficulty it is a logical case to tackle next.

By work of Exoo \cite{Ex90} we know that $R(5, 4.5) \geq 30$ and by unpublished work of Boza we know that $R(5, 4.5) \leq 33$, so it suffices to prove that the set $\cR(5, 4.5, 30)$ is empty.

We briefly describe a graph in $\cR(5,4.5,29)$, using the graph6 format:
{\small
\begin{verbatim}
\~zLugNJije^bwdeq[dIhYX`kWQtPJ_yReEvmSuuQ^ZPRnr\W^cLkzWbNUzIL]vK??~~{
\end{verbatim}
}
We suspect it is the only one, but proving it would require at least one order of magnitude more calculations so we did not do it. We have also found 2\,366\,075 graphs in $\cR(5, 4.5, 28)$ and we suspect there are at least a few more.

Our proof uses a combination of linear programming and computer-assisted brute force calculations, with the key input being a census of $\cR(4, 4.5)$ and $\cR(5, 3.5)$.

We did all the necessary calculations twice, using slightly different methods. The first time our method was quite complicated and used a number of different programs we wrote in C, including a special purpose SAT solver. The second time we kept things as simple as possible, and mainly used a single program written in Rust. A small number of cases did not finish in a reasonable time, and for those we used the open source SAT solver Kissat \cite{kissat} written by Biere. We also used the \emph{Nauty} software package \cite{McPi14} extensively.

The second time around was also faster, with a total running time for proving Theorem \ref{t:main2} of a few months of CPU time. Because the second method was less complicated and has less potential for bugs, we describe this method in detail and briefly outline the first method at the very end of the paper.

\subsection{Acknowledgements}
The author would like to thank Brendan McKay for some help earlier in this project.

\section{A census of $\cR(4,4.5)$ and $\cR(5,3.5)$}
For each of $\cR(4,4.5)$ and $\cR(5,3.5)$ the number of graphs is small enough that we can find all of them by using a one vertex extender. This took a few hours of CPU time in total. See Table \ref{fig:census}.

\begin{table}
\begin{center}
  \begin{tabular}{ l | c | c | c | c }
    $n$ & $|\cR(4, 4.5, n)|$ & edges & $|\cR(5,3.5)|$ & edges \\ \hline \hline
    1 & 1               & [0,0]   & 1         & [0,0] \\ \hline
    2 & 2               & [0,1]   & 2         & [0,1] \\ \hline
    3 & 4               & [0,3]   & 4         & [0,3] \\ \hline
    4 & 10              & [0,5]   & 9         & [2,6] \\ \hline
    5 & 27              & [2,8]   & 21        & [3,9] \\ \hline
    6 & 104             & [3,12]  & 63        & [6,13] \\ \hline
    7 & 509             & [5,16]  & 210       & [9,18] \\ \hline
    8 & 3\,690          & [7,21]  & 897       & [12,24] \\ \hline
    9 & 37\,377         & [9,27]  & 4\,463    & [18,30] \\ \hline
    10 & 480\,802       & [14,33] & 23\,577   & [24,37] \\ \hline
    11 & 6\,303\,192    & [19,40] & 97\,796   & [30,45] \\ \hline
    12 & 60\,750\,166   & [24,48] & 180\,813  & [36,54] \\ \hline
    13 & 265\,108\,612  & [30,52] & 46\,510   & [44,57] \\ \hline
    14 & 323\,189\,897  & [35,60] & 856       & [52,60] \\ \hline
    15 & 64\,354\,692   & [42,64] & 13        & [60,65] \\ \hline
    16 & 1\,494\,223    & [48,72] &           & \\ \hline
    17 & 3\,033         & [61,72] &           & \\ \hline
    18 & 6              & [72,77] &           & \\ \hline \hline
  \end{tabular}
\end{center}
\caption{A census of $\cR(4, 4.5)$ and $\cR(5,3.5)$} \label{fig:census}
\end{table}

An algorithm that extends an $\cR(s, t)$-graph by a single vertex in all possible ways while staying within $\cR(s,t)$ is straightforward, using a variation of the algoritm given in \cite[Section 4]{McRa95}. For example, to extend $G \in \cR(4, 4.5, n)$ by a single vertex we consider the following subsets of $V(G)$:
\begin{enumerate}
 \item A triangle $K_3$.
 \item An independent $3$-set $\overline{K}_3$ that extends to an independent $4$-set.
 \item A $\overline{K}_{3.5}$.
\end{enumerate}
In the case of $\cR(4, 4.5)$, the neighbours of the extra vertex cannot cover a $K_3$, and the complement cannot cover an independent $3$-set that extends to an independent $4$-set or a $\overline{K}_{3.5}$. The case of $\cR(5, 3.5)$ is similar.

\section{Some linear programming} \label{s:LP}
Given a graph $F$ and a vertex $v \in V(F)$, let $F_v^+$ denote the induced subgraph on the vertices adjacent to $v$ and let $F_v^-$ denote the induced subgraph on the vertices not adjacent to $v$ (not including $v$). If $F \in R(s,t,n)$ and $v$ has degree $d$ then $F_v^+ \in R(s-1,t,d)$ and $F_v^- \in R(s,t-1,n-d-1)$. Recall, e.g.\ from \cite{McRa97}, that if $F$ has $n$ vertices we have an equation
\[
 \sum_{v \in V(F)} 2e(F_v^-) = \sum_{v \in V(F)} v(F_v^+)(n-2v(F_v^+)) + 2e(F_v^+).
\]
Here $v(X)$ denotes the number of vertices in $X$ and $e(X)$ denotes the number of edges in $X$. We will refer to this equation as the \emph{edge equation}. Here $v(X)$ denotes the number of vertices in $X$ and $e(X)$ denotes the number of edges in $X$. In particular there must be some vertex $v \in V(F)$ for which
\begin{equation} \label{eq:edge}
 2e(F_v^-) \leq v(F_v^+)(n-2v(F_v^+)) + 2e(F_v^+).
\end{equation}
We will use this equation repeatedly. If $F \in \cR(5, 4.5, 30)$ then any vertex $v \in V(F)$ must have degree $d(v) \in \{14, 15, 16, 17, 18\}$, and the edge equation says that we must have some vertex $v \in V(F)$ with
\begin{eqnarray*}
 e(F_v^-) & \leq & e(F_v^+) + 14 \quad \textnormal{if $v(F_v^+) = 14$} \\
 e(F_v^-) & \leq & e(F_v^+) \phantom{ {} + 14} \quad \textnormal{if $v(F_v^+) = 15$} \\
 e(F_v^-) & \leq & e(F_v^+) - 16 \quad \textnormal{if $v(F_v^+) = 16$} \\
 e(F_v^-) & \leq & e(F_v^+) - 34 \quad \textnormal{if $v(F_v^+) = 17$} \\
 e(F_v^-) & \leq & e(F_v^+) - 54 \quad \textnormal{if $v(F_v^+) = 18$} \\
\end{eqnarray*}
Here $F_v^+$ is in $\cR(4, 4.5)$ and $F_v^-$ is in $\cR(5, 3.5)$, and we will use this to reduce the number of gluing operations needed.

This is helpful because a typical graph in $\cR(4, 4.5)$ is less dense while a typical graph in $\cR(5, 3.5)$ is more dense, so most pairs $(G, H)$ with $G \in \cR(5, 3.5, 30 - 1 - d)$ and $H \in \cR(4, 4.5, d)$ will not satisfy the above inequality.

\section{Gluing} \label{s:gluing}
We used a method similar to the one described in \cite[Section 3]{McRa95}. This whole section is a summary of that section, with a few modifications to the current situation. Suppose that $G \in \cR(5, 3.5)$ and $H \in \cR(4, 4.5)$. We want to find all $\cR(5, 4.5)$-graphs $F$ for which some vertex $x \in V(F)$ has $F_x^- = G$ and $F_x^+ = H$.

Define a \emph{feasible cone} to be a subset $C$ of $V(G)$ which covers no clique of size 4, and whose complement covers no independent set of size 3. The point is that the neighbourhood $N_F(v, V(G))$ of $v$ in $G$ must be a feasible cone for every $v \in V(H)$. The condition that the complement covers no independent set of size 3 comes from the fact that if the complement did cover an independent set $\{a, b, c\}$ of size 3 then $\{a, b, c, v, x\}$ would be a $\overline{K }_{4.5}$.

Now the problem becomes finding all ways to choose feasible cones $C_0, C_1,\ldots$, one for each vertex of $H$, such that no $K_5$ or $\overline{K}_{4.5}$ appear in $F$. We need to consider all the ways this can happen. Let $\overline{C}$ denote the complement of $C$ in $V(G)$.
\begin{eqnarray*}
 & K_2: & \textnormal{Two adjacent vertices $v, w \in V(H)$ have $C_v \cap C_w$ covering some triangle} \\ & & \textnormal{in $G$.} \\
 & K_3: & \textnormal{Three adjacent vertices $u, v, w \in V(H)$ have $C_u \cap C_v \cap C_w$ covering some} \\ & & \textnormal{edge in $G$.} \\
 & E_2: & \textnormal{Two non-adjacent vertices $v, w \in V(H)$ have $\overline{C}_v \cap \overline{C}_w$ covering some $\overline{K}_{2.5}$} \\ & &  \textnormal{in $G$.} \\
 & E_3: & \textnormal{Three non-adjacent vertices $u, v, w \in V(H)$ have $\overline{C}_u \cap \overline{C}_v \cap \overline{C}_w$ covering} \\ & &  \textnormal{some edge in $G$.} \\
 & E_{2.5}: & \textnormal{Three vertices $u, v, w \in V(H)$ with a single edge between them have} \\ & & \textnormal{$\overline{C}_u \cap \overline{C}_v \cap \overline{C}_w$ covering some independent $2$-set in $G$.} \\
 & E_{3.5}: & \textnormal{Four vertices $t, u, v, w \in V(H)$ with a single edge between them have} \\ & & \textnormal{$\overline{C}_t \cap \overline{C}_u \cap \overline{C}_v \cap \overline{C}_w$ non-empty.} \\
 & E_3' & \textnormal{Three non-adjacent vertices $u, v, w \in V(H)$ have $\overline{C}_u \cap \overline{C}_v$ covering some} \\ & & \textnormal{independent $2$-set $\{a, b\}$ in $G$ and $\overline{C}_w$ covers either $a$ or $b$.} \\
 & E_4' & \textnormal{Four non-adjacent vertices $t, u, v, w \in V(H)$ have $\overline{C}_t \cap \overline{C}_u \cap \overline{C}_v$} \\ & & \textnormal{non-empty.}
\end{eqnarray*}

In the last two conditions $E_3'$ and $E_4'$, the vertex $w$ can be any of the vertices in the independent $3$-set or $4$-set.

Note that we do not need a case for two non-adjacent vertices $v, w \in V(H)$ with $\overline{C}_v \cap \overline{C}_w$ covering a $\overline{K}_3$ in $G$.

To reduce the number of cases we consider an \emph{interval}, or \emph{interval of feasible cones}, to be a set of feasible cones of the form $\{X \quad | \quad B \subset X \subset T \}$ for some feasible cones $B \subset T$. This interval $[B, T]$ contains $2^{|T|-|B|}$ feasible cones.

We now have 8 collapsing rules, analogous to the 4 collapsing rules (a) - (d) described in \cite[p.\ 5]{McRa95} with each corresponding to one of the 8 ways a $K_5$ or $\overline{K}_{4.5}$ can appear. We describe one of them in detail, leaving the rest to the reader:
\newline

\noindent
($E_3'$ collapsing rule) We define a function $G_2 : 2^{V(G)} \to 2^{V(G)}$ with
\[
 G_2(X) = \{b \in V(G) \quad | \quad \{a, b\} \textnormal{ is an independent $2$-set for some $a \in X$}.\}
\]
Suppose $\{u, v, w\}$ is an independent $3$-set in $H$. Let $Y = G_2(\overline{T}_u \cap \overline{T}_v \cap \overline{T}_w)$.

\textbf{if} $Y \cap \overline{T}_u \cap \overline{T}_v \neq \emptyset$ \textbf{then} FAIL

\textbf{if} $Y \cap \overline{T}_u \cap \overline{T}_w \neq \emptyset$ \textbf{then} FAIL

\textbf{if} $Y \cap \overline{T}_v \cap \overline{T}_w \neq \emptyset$ \textbf{then} FAIL

\textbf{else}

\hspace{20pt} $B_u := B_u \cup (Y \cap \overline{T}_v) \cup (Y \cap \overline{T}_w)$

\hspace{20pt} $B_v := B_v \cup (Y \cap \overline{T}_u) \cup (Y \cap \overline{T}_w)$

\hspace{20pt} $B_w := B_w \cup (Y \cap \overline{T}_u) \cup (Y \cap \overline{T}_v)$

The reason for the FAIL is that if, for example, $b \in Y \cap \overline{T}_u \cap \overline{T}_v$ then there exists a vertex $a \in V(G)$ not adjacent to $b$ or any of $u$, $v$ and $w$, and the only possible edge between the vertices $\{a, b, u, v, w\}$ is the single edge $(b, w)$. Hence we are guaranteed a $\overline{K}_{4.5}$ or $\bar{K}_5$.

Similarly, if $b \in Y \cap \overline{T}_v$ then there exists a vertex $a \in V(G)$ not adjacent to $b$ or any of $u$, $v$ and $w$, and the only two possible edges between the vertices $\{a, b, u, v, w\}$ are $(b, u)$ and $(b, w)$. Hence those must both be present to avoid a $\overline{K}_{4.5}$.
\newline

All the remaining collapsing rules are similar.

The analogues of all the remaining results in \cite[Section 3]{McRa95} hold in our setting, and each collapsing rule can be carried out using a small number of bitwise operations and a precomputed table like the function $G_2$ described above. We used the same kind of pair of superimposed trees, called a \emph{double tree} in loc.\ cit., to process many $H$ at the same time.

We can also reverse the roles of $G$ and $H$, and this leads a very similar algorithm with similar collapsing rules.

\section{Encoding a gluing problem as a SAT problem} \label{s:SATsolver}
The other basic approach is to encode the problem as a SAT problem. Suppose once again that we have $G \in \cR(5, 3.5)$ and $H \in \cR(4, 4.5)$, and that we want to find all $\cR(5, 4.5)$-graphs $F$ for which some vertex $x \in V(F)$ has $F_x^- = G$ and $F_x^+ = H$. We can encode this as a SAT problem with one variable for each possible edge between a vertex in $G$ and a vertex in $H$, and one clause for each way to make a $K_5$ or a $\overline{K}_{4.5}$.

One can of course construct a single SAT problem that encodes multiple gluing problems, but because we covered most of the cases using the method described in Section \ref{s:gluing} this was not necessary so to keep our method as simple as possible we did not do it.

\section{Detailed graph counts}
Suppose $G \in \cR(5, 4.5, 30)$. Then every vertex of $G$ must have degree in the set $\{14, 15, 16, 17, 18\}$. If a vertex $v$ has degree $d_1$ then $G_v^+ \in \cR(4, 4.5, d_1)$ and $G_v^- \in \cR(5, 3.5, d_2)$ with $d_2 = 30 - d_1 - 1$, so
\[
 (d_1, d_2) \in \{(14, 15), (15, 14), (16, 13), (17, 12), (18, 11)\}.
\]
We go through each pair $(d_1, d_2)$ in turn:

\subsection{$(d_1, d_2) = (14, 15)$}
For $(d_1, d_2) = (14, 15)$ the counts are as in Table \ref{t:1415}:
\begin{table}[h]
\begin{center}
  \begin{tabular}{ l | c | c | c | }
    e & $|\cR(5, 3.5, 15, e)|$ & e & $|\cR(4, 4.5, 14, e)|$ \\ \hline \hline
    60 & 1    & 46 & 56\,042\,064 \\ \hline
    61 & 2    & 47 & 29\,676\,081 \\ \hline
    62 & 4    & 48 & 11\,623\,714 \\ \hline
    63 & 5    & 49 & 3\,558\,725 \\ \hline
    64 & 0    & 50 & 925\,641 \\ \hline
    65 & 1    & 51 & 223\,034 \\ \hline
       &      & 52 & 52\,386 \\ \hline
       &      & 53 & 11\,909 \\ \hline
       &      & 54 & 2\,663 \\ \hline
       &      & 55 & 550 \\ \hline
       &      & 56 & 127 \\ \hline
       &      & 57 & 23 \\ \hline
       &      & 58 & 7 \\ \hline
       &      & 59 & 1 \\ \hline
       &      & 60 & 1 \\ \hline
\end{tabular}
\end{center}
\caption{A partial census of $\cR(5, 3.5, 15)$ and $\cR(4, 4.5, 14)$}
\label{t:1415}
\end{table}
Here we have only included $|\cR(4, 4.5, 14, e)|$ for those $e$ necessary according to Equation \ref{eq:edge}.

Without any simplifications this yields about 280 million gluing operations. But using the method described in Section \ref{s:gluing} above with a single graph $G \in \cR(5, 3.5, 15)$ and all the graphs in $\cR(4, 4.5, 14)$ organised into a tree was quite fast so we glued all of them. Here the memory consumption of our program became a problem, so we processed the 322 million graphs in $\cR(4, 4.5, 14)$ in batches of about 27 million at a time. (It is almost certainly possible to make our program more memory efficient, but it was only a problem in this one case and it did not seem worth the trouble. Such a double tree with 27 million graphs used about 13GB of memory, which was the most our computer could handle.)

\FloatBarrier

\subsection{$(d_1, d_2) = (15, 14)$}
For $(d_1, d_2) = (15, 14)$ we have a similar census, given in Table \ref{t:1514}:
\begin{table}[h]
\begin{center}
  \begin{tabular}{ l | c | c | c | }
    e & $|\cR(5, 3.5, 14, e)|$ & e & $|\cR(4, 4.5, 15, e)|$ \\ \hline \hline
    52 & 2    & 52 & 15\,126\,241 \\ \hline
    53 & 14   & 53 & 10\,096\,214 \\ \hline
    54 & 55   & 54 & 4\,670\,773 \\ \hline
    55 & 133  & 55 & 1\,540\,615 \\ \hline
    56 & 191  & 56 & 384\,517 \\ \hline
    57 & 198  & 57 & 82\,242 \\ \hline
    58 & 132  & 58 & 18\,260 \\ \hline
    59 & 78   & 59 & 4\,523 \\ \hline
    60 & 35   & 60 & 1\,218 \\ \hline
    61 & 12   & 61 & 288 \\ \hline
    62 & 4    & 62 & 66 \\ \hline
    63 & 2    & 63 & 10 \\ \hline
       &      & 64 & 2 \\ \hline
\end{tabular}
\end{center}
\caption{A partial census of $\cR(5, 3.5, 14)$ and $\cR(4, 4.5, 15)$}
\label{t:1514}
\end{table}

Without any simplifications this yields just over one billion gluing operations. Once again we used the method described in Section \ref{s:gluing} to glue each $G \in \cR(5, 3.5, 14)$ to a tree of graphs in $\cR(4, 4.5, 15)$. The running time is quite sensitive to the number of graphs of order 9 (and to some extent 10) in the tree, and also very sensitive to how dense those graphs are. Hence we found it convenient to use a tree that contained almost all the necessary $\cR(4, 4.5, 15)$-graphs while leaving out a small number. In the end we left out 1\,335 graphs in $\cR(4, 4.5, 15, e \geq 52)$.

With those graphs, we used the method described in Section \ref{s:gluing} with the roles of $G$ and $H$ reversed.

Most of these remaining cases were quite fast, but a few were very slow and we did not run the slowest cases to completion.

Finally, we processed the remaining pairs of graphs with the SAT solver Kissat \cite{kissat} using the method described in Section \ref{s:SATsolver}.

\subsection{$(d_1, d_2) = (16, 13)$}
For $(d_1, d_2) = (16, 13)$ we have the counts given in Table \ref{t:1613}:
\begin{table}[h]
\begin{center}
  \begin{tabular}{ l | c | c | c | }
    e & $|\cR(5, 3.5, 13, e)|$ & e & $|\cR(4, 4.5, 16, e)|$ \\ \hline \hline
    44 & 7    & 60 & 275\,651 \\ \hline
    45 & 69   & 61 & 131\,157 \\ \hline
    46 & 518  & 62 & 41\,670 \\ \hline
    47 & 2\,271 & 63 & 9\,171 \\ \hline
    48 & 6\,246 & 64 & 1\,494 \\ \hline
    49 & 10\,903 & 65 & 225 \\ \hline
    50 & 12\,235 & 66 & 73 \\ \hline
    51 & 8\,774 & 67 & 34 \\ \hline
    52 & 4\,008 & 68 & 29 \\ \hline
    53 & 1\,191 & 69 & 9 \\ \hline
    54 & 245  & 70 & 5 \\ \hline
    55 & 36   & 71 & 1 \\ \hline
    56 & 6    & 72 & 1 \\ \hline
    57 & 1    &  &  \\ \hline
\end{tabular}
\end{center}
\caption{A partial census of $\cR(5, 3.5, 13)$ and $\cR(4, 4.5, 16)$}
\label{t:1613}
\end{table}

Without any simplifications this yields about 87 million gluing operations. Now we used a version of the algorithm described in Section \ref{s:gluing} with the roles of $G$ and $H$ reversed, with a dense subgraph of $G$ in front. We thought this might be quite slow because we have 459\,520 graphs in $\cR(4, 4.5, 16)$ to consider but the graphs in $\cR(4, 4.5, 16, e)$ ran faster for smaller $e$, in part because we needed to consider fewer graphs in $\cR(5, 3.5, 13)$, so that turned out not to be a problem. Once again we had a few leftover graphs that we processed with a SAT solver.

\subsection{$(d_1, d_2) = (17, 12)$}
For $(d_1, d_2) = (17, 12)$ we have consider Table \ref{t:1712}.
\begin{table}[h]
\begin{center}
  \begin{tabular}{ l | c | c | c | }
    e & $|\cR(5, 3.5, 12, e \leq 38)|$ & e & $|\cR(4, 4.5, 17, e)|$ \\ \hline \hline
    36 & 4   & 70 & 96 \\ \hline
    37 & 29  & 71 & 10 \\ \hline
    38 & 365 & 72 & 1 \\ \hline
\end{tabular}
\end{center}
\caption{A partial census of $\cR(5, 3.5, 12)$ and $\cR(4, 4.5, 17)$}
\label{t:1712}
\end{table}

Because this was quite fast anyway we glued all of $\cR(5, 3.5, 12)$ to all of $\cR(4, 4.5, 17)$ using the method in Section \ref{s:gluing} with the roles of $G$ and $H$ reversed.

\subsection{$(d_1, d_2) = (18, 12)$}
Finally, for $(d_1, d_2) = (18, 12)$ we have $e(\cR(5, 3.5, 11)) \geq 30$ and $e(\cR(4, 4.5, 18)) \leq 77$, so for this case there is nothing to do. Because it was very fast anyway, we glued all of $\cR(5, 3.5, 11)$ to all of $\cR(4, 4.5, 18)$ using the method in Section \ref{s:gluing} with the roles of $G$ and $H$ reversed.

None of our gluing operations yielded any outputs, and all of the SAT instances returned UNSAT, so that completes the proof of Theorem \ref{t:main2}.

\section{The first method}
We briefly describe the more complicated method we used the first time. We considered various \emph{intervals of graphs}, where an interval $[B, T]$ of graphs with a fixed vertex set consists of all graphs $X$ with $E(B) \subset E(X) \subset E(T)$.

We used a simple SAT solver written in C which, depending on the order of the possible solution graphs, output either all solutions or a set of intervals of possible solutions.

\subsection{$(d_1, d_2) = (14, 15)$}
We covered $\cR(5, 3.5, 15)$ by $6$ intervals in $\cR(5, 3.5, 15)$ and we covered $\cR(4, 4.5, 14)$ by $819$ intervals in $\cR(4, 4.5, 11)$. This yielded a total of $1\,904$ gluing operations (as not all the intervals had to be glued to each other).

This produced $2\,450\,116$ graphs of order $27$. These extended to $197\,411$ graphs of order $28$, then none of order $29$.

\subsection{$(d_1, d_2) = (15, 14)$}
We covered $\cR(5, 3.5, 14)$ by $140$ intervals in $\cR(5, 3.5, 14)$. Then we covered each set $\cR(4, 4.5, 15, e \geq e_0)$ by some intervals in $\cR(4, 4.5, 10)$ and some leftover graphs in $\cR(4, 4.5, 15)$. For example, for $\cR(4, 4.5, 15, e \geq 52)$ we used $53$ intervals in $\cR(4, 4.5, 10)$ and $165$ graphs in $\cR(4, 4.5, 15)$.

This produced $112\,159\,105$ graphs of order $27$. These extended to $1\,252\,566$ graphs of order $28$, then none of order $29$.

\subsection{$(d_1, d_2) = (16, 13)$}
We covered $\cR(5, 3.5, 13)$ by $6$ intervals in $\cR(5, 3.5, 8)$, $34$ intervals in $\cR(5, 3.5, 11)$, $29$ intervals in $\cR(5, 3.5, 12)$, and $420$ intervals in $\cR(5, 3.5, 13)$. Then, depending on the type of interval in $\cR(5, 3.5)$, we used different types of intervals in $\cR(4, 4.5)$ ranging from $\cR(4, 4.5, 11)$ to $\cR(4, 4.5, 16)$.

Moreover, we left out some gluing operations for $\cR(4, 4.5, 16, e)$ with $\cR(3, 5.5, 13, e-16)$. That left $36879$ graphs in $\cR(4, 4.5, 16)$, all with minimum degree at most $7$, and we glued these along an edge. By that we mean that we considered the union of two $\cR(4, 4.5, 16)$-graphs $G_1$ and $G_2$ with distinguished vertices $a \in V(G_1)$ and $b \in V(G_2)$ along $(G_1)_a^+ \cong (G_2)_b^+$. This was similar to the metod described in \cite{AnMc18, AnMc25}.

\subsection{$(d_1, d_2) = (17, 12)$}
We covered $\cR(5, 3.5, 12, e \leq 38)$ by $7$ intervals in $\cR(5, 3.5, 9)$. This reduced the calculation to $229$ gluing operations. This produced a total of $82\,851$ graphs of order 27. These extended to $618$ graphs of order 28, then none of order 29.

\section{Concluding remarks}
It is probably feasible to compute the set $\cR(5, 4.5, 29)$ of maximal Ramsey graphs, but it would take significantly longer so we did not do it.



\end{document}